\documentclass[12pt]{article}
\newcommand{\nc}{\newcommand}
\nc{\dfrac}{\displaystyle \frac } 
\newcommand{\beq}{\begin{equation}}
\newcommand{\eeq}{\end{equation}}
\newcommand{\br}{\begin{array}}
\newcommand{\er}{\end{array}}

\input{tcilatex}

\begin{document}

\begin{center}
{\Large \textbf{Quantum Group Covariance and the Braided Structure of
Deformed Oscillators }}
\end{center}

\vspace{5mm}

\begin{center}
A. Yildiz

\noindent Feza G\"ursey Institute, P.O. Box 6, 81220, \c{C}engelk\"{o}y,
Istanbul, Turkey \footnote{%
E-mail : yildiz@gursey.gov.tr}. \\[0pt]
\end{center}

\vspace{5mm}

\begin{center}
\textbf{Abstract}
\end{center}

The connection between braided Hopf algebra structure and quantum group
covariance of the deformed oscillators is constructed explicitly. In this
context we provide deformations of the Hopf algebra of functions on $SU(1,1).
$ Quantum subgroups and their representations are also discussed.

\noindent

\vspace{1cm} \noindent

\section{Introduction}

The covariance of the oscillator algebras attracted a lot of attention and
is discussed in different contexts$^{1}$. The covariance of an algebra under
the action of a noncommutative algebra deforms the notion of defining
identical copies in the transformed algebra and this leads to the
deformation of the usual tensor product namely braided tensor product. The
Hopf algebra axioms are replaced by the braided Hopf algebra axioms$^{2}$.
Hence braided group theory (self contained reviews can be found in Ref.3)
unifies the notions of symmetry and statistics. Recently, we found the
general braided Hopf algebra solutions of the generalized oscillators$^{4}$.
In this work we show that some of these solutions are connected with the
quantum group covariance and we find the $\mathbf{R}$-matrices controlling
the braiding structure and the quantum group. We also discuss the
representations of quantum subgroups.

\section{The generalized oscillator, its covariance and braided Hopf
structure}

Suppose that the generalized oscillator algebra

\begin{eqnarray}
aa^{\ast }-Q_{1}a^{\ast }a &=&q^{2N}  \nonumber  \label{eq:osc1} \\
aq^{N} &=&qq^{N}a \\
q^{N}a^{\ast } &=&qa^{\ast }q^{N}  \nonumber
\end{eqnarray}

\noindent is covariant under the transformation
\begin{eqnarray}  \label{eq:tfm1}
(a)^{\prime }&=&aK_{1}+q^{N}K_{2}+a^{\ast }K_{3}  \nonumber \\
(a^{\ast })^{\prime }&=&a^{\ast }K_{1}^{\ast }+q^{N}K_{2}^{\ast
}+aK_{3}^{\ast } \\
(q^{N})^{\prime }&=&aL_{1}+q^{N}L_{2}+a^{\ast }L_{1}^{\ast }  \nonumber
\end{eqnarray}

\noindent where the deformation parameters $Q_{1}$and $q$ are positive and
with the $\ast $-structure $(a^{\ast })^{\ast }=a$ and $(q^{N})^{\ast }=q^{N}
$. The elements $K_{1},K_{1}^{\ast },K_{2},K_{2}^{\ast },K_{3},K_{3}^{\ast
},L_{1},L_{1}^{\ast }$ and $L_{2}$ ($L_{2}^{\ast }=L_{2}$) generate some
algebra. Our aim is to find that algebra if the transformation is a quantum
group transformation. We write the above transformation as a covector
transformation

\begin{equation}  \label{eq:tfm2}
\mathbf{x}^{\prime }\mathbf{=xt}
\end{equation}

\noindent where

\begin{equation}  \label{eq:qm1}
\mathbf{x}=(
\begin{array}{ccc}
a & a^{\ast } & q^{N}
\end{array}
),\ \ \ \mathnormal{and}\ \ \ \mathbf{t}=\left(
\begin{array}{ccc}
K_{1} & K_{3}^{\ast } & L_{1} \\
K_{3} & K_{1}^{\ast } & L_{1}^{\ast } \\
K_{2} & K_{2}^{\ast } & L_{2}
\end{array}
\right) .
\end{equation}

\noindent The matrix $\mathbf{t}$ is a quantum matrix satisfying
\begin{equation}  \label{eq:qg1}
\mathbf{Rt_{1}{t}_{2}=t_{2}{t}_{1}R}
\end{equation}

\noindent and $\mathbf{R}$ satisfies QYBE

\begin{equation}  \label{eq:QYBE}
\mathbf{R_{12}R_{13}R_{23}}=\mathbf{R_{23}R_{13} R_{12}}.
\end{equation}

\noindent To find the $\mathbf{R}$- matrix and hence the quantum group, we
write the oscillator algebra as a covector algebra
\begin{equation}  \label{eq:covalg}
\mathbf{x_1}\mathbf{x_2}=\mathbf{x_2}\mathbf{x_{1}R}
\end{equation}

\noindent where

\begin{equation}
\mathbf{x_1}\mathbf{x_2}=\left(
\begin{array}{ccccccccc}
a^{2} & aa^{\ast } & aq^{N} & a^{\ast }a & (a^{\ast })^{2} & a^{\ast }q^{N}
& q^{N}a & q^{N}a^{\ast } & q^{2N}
\end{array}
\right)
\end{equation}
\begin{equation}
\mathbf{x}_{2}\mathbf{x_{1}}=\left(
\begin{array}{ccccccccc}
a^{2} & a^{\ast }a & q^{N}a & aa^{\ast } & (a^{\ast })^{2} & q^{N}a^{\ast }
& aq^{N} & a^{\ast }q^{N} & q^{2N}
\end{array}
\right)
\end{equation}

\noindent and the general form of the $\mathbf{R}$-matrix

\begin{equation}  \label{eq:rm1}
\mathbf{R}=\left(
\begin{array}{ccccccccc}
1 & 0 & 0 & 0 & 0 & 0 & 0 & 0 & 0 \\
0 & A_{1} & 0 & A_{6} & 0 & 0 & 0 & 0 & A_{15} \\
0 & 0 & A_{4} & 0 & 0 & 0 & A_{11} & 0 & 0 \\
0 & A_{2} & 0 & A_{7} & 0 & 0 & 0 & 0 & A_{16} \\
0 & 0 & 0 & 0 & 1 & 0 & 0 & 0 & 0 \\
0 & 0 & 0 & 0 & 0 & A_{9} & 0 & A_{13} & 0 \\
0 & 0 & A_{5} & 0 & 0 & 0 & A_{12} & 0 & 0 \\
0 & 0 & 0 & 0 & 0 & A_{10} & 0 & A_{14} & 0 \\
0 & A_{3} & 0 & A_{8} & 0 & 0 & 0 & 0 & A_{17}
\end{array}
\right).
\end{equation}

\noindent The constants ($A_{i}$) appearing in the $\mathbf{R}$-matrix is to
be determined from the consistency of (\ref{eq:covalg}) with the oscillator
relations (\ref{eq:osc1}) and from (\ref{eq:QYBE}). The covariance of a
covector algebra under the action of a quantum group induces a braided Hopf
algebra structure whose axioms are collectively given by

\begin{eqnarray}
m\circ (id\otimes m) &=&m\circ (m\otimes id)  \nonumber  \label{eq:bhopf} \\
m\circ (id\otimes \eta ) &=&m\circ (\eta \otimes id)=id  \nonumber \\
(id\otimes \Delta )\circ \Delta &=&(\Delta \otimes id)\circ \Delta  \nonumber
\\
(\epsilon \otimes id)\circ \Delta &=&(id\otimes \epsilon )\circ \Delta =id
\nonumber \\
m\circ (id\otimes S)\circ \Delta &=&m\circ (S\otimes id)\circ \Delta =\eta
\circ \epsilon  \nonumber \\
\psi \circ (m\otimes id) &=&(id\otimes m)\circ (\psi \otimes id)\circ
(id\otimes \psi )  \nonumber \\
\psi \circ (id\otimes m) &=&(m\otimes id)\circ (id\otimes \psi )\circ (\psi
\otimes id)  \nonumber \\
(id\otimes \Delta )\circ \psi &=&(\psi \otimes id)\circ (id\otimes \psi
)\circ (\Delta \otimes id)  \nonumber \\
(\Delta \otimes id)\circ \psi &=&(id\otimes \psi )(\psi \otimes id)\circ
(id\otimes \Delta )  \nonumber \\
\Delta \circ m &=&(m\otimes m)(id\otimes \psi \otimes id)\circ (\Delta
\otimes \Delta ) \\
S\circ m &=&m\circ \psi \circ (S\otimes S)  \nonumber \\
\Delta \circ S &=&(S\otimes S)\circ \psi \circ \Delta  \nonumber \\
\epsilon \circ m &=&\epsilon \otimes \epsilon  \nonumber \\
(\psi \otimes id)\circ (id\otimes \psi )\circ (\psi \otimes id)
&=&(id\otimes \psi )\circ (\psi \otimes id)\circ (id\otimes \psi ).
\nonumber
\end{eqnarray}

\noindent The $\ast $-structure for a braided algebra $B$ is different from
the non-braided one such that

\begin{eqnarray}
\Delta \circ \ast &=&\pi \circ (\ast \otimes \ast )\circ \Delta  \nonumber \\
S\circ \ast &=&\ast \circ S \\
(a\otimes b)^{\ast } &=&b^{\ast }\otimes a^{\ast },\ \forall a,b\in B.
\nonumber
\end{eqnarray}
\bigskip

\noindent The braided covector algebra has a braided Hopf algebra structure
\begin{equation}
\Delta (\mathbf{x})=\mathbf{x}\otimes 1+1\otimes \mathbf{x},\ \ \epsilon (%
\mathbf{x})=0,\ \ S(\mathbf{x})=-\mathbf{x}
\end{equation}

\noindent with the braiding relations

\begin{equation}
\psi (\mathbf{x}_{1}\otimes \mathbf{x}_{2})=\mathbf{x}_{2}\otimes \mathbf{x}%
_{1}\mathbf{R}^{\prime },\quad \mathnormal{i.e.},\quad \psi (x_{i}\otimes
x_{j})=x_{b}\otimes x_{a}R_{ij}^{\prime ab}.
\end{equation}

\noindent The matrix $\mathbf{R}^{\prime }$ which controls the braiding
relations should satisfy the following conditions

\begin{eqnarray}
\mathbf{R_{12}^{\prime }R_{13}^{\prime }R_{23}^{\prime }}&=& \mathbf{%
R_{23}^{\prime }R_{13}^{\prime }R_{12}^{\prime }}  \nonumber \\
\mathbf{R_{12}^{\prime }R_{13}^{\prime }R_{23}}&=&\mathbf{R
_{23}R_{13}^{\prime }R_{12}^{\prime }}  \nonumber \\
\mathbf{R_{12}R_{13}^{\prime }R_{23}^{\prime }}&=&\mathbf{R _{23}^{\prime
}R_{13}^{\prime }R_{12}} \\
\mathbf{(PR^{\prime }+1)(PR-1)}&=&\mathbf{0}  \nonumber \\
\mathbf{R_{21}^{\prime }R}&=&\mathbf{R_{21}R}  \nonumber
\end{eqnarray}

\noindent where $\mathbf{P}$ is the permutation matrix. Hence the problem of
finding the quantum group leaving the generalized oscillator algebra
covariant and the braidings induced by the quantum group is reduced to
finding the matrices $\mathbf{R}$ and $\mathbf{R}^{\prime }$.

The general form of the matrix $\mathbf{R}^{\prime }$ can be written as
\begin{equation}  \label{eq:rmprime1}
\mathbf{R}^{\prime }=\left(
\begin{array}{ccccccccc}
C_{1} & 0 & 0 & 0 & 0 & 0 & 0 & 0 & 0 \\
0 & C_{2} & 0 & C_{7} & 0 & 0 & 0 & 0 & C_{12} \\
0 & 0 & C_{5} & 0 & 0 & 0 & C_{11} & 0 & 0 \\
0 & C_{3} & 0 & C_{8} & 0 & 0 & 0 & 0 & C_{13} \\
0 & 0 & 0 & 0 & C_{1} & 0 & 0 & 0 & 0 \\
0 & 0 & 0 & 0 & 0 & C_{10} & 0 & C_{6} & 0 \\
0 & 0 & C_{6} & 0 & 0 & 0 & C_{10} & 0 & 0 \\
0 & 0 & 0 & 0 & 0 & C_{11} & 0 & C_{5} & 0 \\
0 & C_{4} & 0 & C_{9} & 0 & 0 & 0 & 0 & C_{14}
\end{array}
\right)
\end{equation}

\noindent which gives the general form of the braiding relations.

For the three deformation parameters $Q_{1}$ and $q$ free, it turns out that
there is a unique solution for the matrices $\mathbf{R}$ and $\mathbf{R}%
^{\prime }$ namely

\begin{equation}
R=\left(
\begin{array}{ccccccccc}
1 & 0 & 0 & 0 & 0 & 0 & 0 & 0 & 0 \\
0 & \dfrac{Q_{1}^{2}}{q^{2}} & 0 & 0 & 0 & 0 & 0 & 0 & 0 \\
0 & 0 & \dfrac{Q_{1}}{q} & 0 & 0 & 0 & 0 & 0 & 0 \\
0 & \dfrac{(q^{2}-Q_{1})}{q^{2}} & 0 & \dfrac{1}{Q_{1}} & 0 & 0 & 0 & 0 & 0
\\
0 & 0 & 0 & 0 & 1 & 0 & 0 & 0 & 0 \\
0 & 0 & 0 & 0 & 0 & \dfrac{1}{q} & 0 & \dfrac{(q^{2}-Q_{1})}{q^{2}} & 0 \\
0 & 0 & \dfrac{(q^{2}-Q_{1})}{q^{2}} & 0 & 0 & 0 & \dfrac{1}{q} & 0 & 0 \\
0 & 0 & 0 & 0 & 0 & 0 & 0 & \dfrac{Q_{1}}{q} & 0 \\
0 & \dfrac{Q_{1}}{q^{2}} & 0 & -\dfrac{1}{Q_{1}} & 0 & 0 & 0 & 0 & 1
\end{array}
\right)
\end{equation}

\noindent and

\begin{equation}
\mathbf{R}^{\prime }=\left(
\begin{array}{ccccccccc}
\dfrac{q^{2}}{Q_{1}} & 0 & 0 & 0 & 0 & 0 & 0 & 0 & 0 \\
0 & Q_{1} & 0 & 0 & 0 & 0 & 0 & 0 & 0 \\
0 & 0 & q & 0 & 0 & 0 & 0 & 0 & 0 \\
0 & \dfrac{q^{2}-Q_{1}}{Q_{1}} & 0 & \dfrac{q^{2}}{Q_{1}^{2}} & 0 & 0 & 0 & 0
& 0 \\
0 & 0 & 0 & 0 & \dfrac{q^{2}}{Q_{1}} & 0 & 0 & 0 & 0 \\
0 & 0 & 0 & 0 & 0 & \dfrac{q}{Q_{1}} & 0 & \dfrac{q^{2}-Q_{1}}{Q_{1}} & 0 \\
0 & 0 & \dfrac{q^{2}-Q_{1}}{Q_{1}} & 0 & 0 & 0 & \dfrac{q}{Q_{1}} & 0 & 0 \\
0 & 0 & 0 & 0 & 0 & 0 & 0 & q & 0 \\
0 & 1 & 0 & -\dfrac{q^{2}}{Q_{1}^{2}} & 0 & 0 & 0 & 0 & \dfrac{q^{2}}{Q_{1}}
\end{array}
\right)   \label{eq:rmprime2}
\end{equation}

\noindent Similar to the $R$ matrix of $SU_{q}(2)$ the matrix $R^{\prime }$
is proportional to $R$ $(R^{\prime }=q^{2}Q_{1}^{-1}R)$. The entries of the
quantum matrix (\ref{eq:qm1}) generate the algebra

\begin{eqnarray}
K_{1}K_{1}^{\ast } &=&K_{1}^{\ast }K_{1}+q^{2}Q_{1}^{-2}L_{1}^{\ast
}L_{1}+q^{-2}Q_{1}(q^{2}-Q_{1})K_{3}^{\ast }K_{3}  \nonumber
\label{eq:relations1} \\
K_{1}K_{2} &=&qQ_{1}^{-1}K_{2}K_{1}  \nonumber \\
K_{1}K_{2}^{\ast } &=&q^{-1}Q_{1}K_{2}^{\ast
}K_{1}+qQ_{1}^{-1}L_{2}L_{1}+q^{-2}Q_{1}(q^{2}-Q_{1})K_{3}^{\ast }K_{2}
\nonumber \\
K_{1}K_{3} &=&q^{2}Q_{1}^{-2}K_{3}K_{1}  \nonumber \\
K_{1}K_{3}^{\ast } &=&Q_{1}K_{3}^{\ast }K_{1}+L_{1}^{2}  \nonumber \\
K_{1}L_{1} &=&qL_{1}K_{1}  \nonumber \\
K_{1}L_{1}^{\ast } &=&qQ_{1}^{-1}L_{1}^{\ast
}K_{1}+q^{-1}(q^{2}-Q_{1})L_{1}K_{3}  \nonumber \\
K_{1}L_{2} &=&L_{2}K_{1}+q^{-1}(q^{2}-Q_{1})L_{1}K_{2}  \nonumber \\
K_{2}K_{2}^{\ast } &=&Q_{1}K_{2}^{\ast }K_{2}+q^{-2}Q_{1}^{2}K_{3}^{\ast
}K_{3}-K_{1}^{\ast }K_{1}+L_{2}^{2}  \nonumber \\
K_{2}K_{3} &=&qQ_{1}^{-1}K_{3}K_{2}  \nonumber \\
K_{2}K_{3}^{\ast } &=&q^{-1}Q_{1}^{2}K_{3}^{\ast }K_{2}+L_{2}L_{1}  \nonumber
\\
K_{2}L_{1} &=&Q_{1}L_{1}K_{2} \\
K_{2}L_{1}^{\ast } &=&L_{1}^{\ast }K_{2}+q^{-1}(q^{2}-Q_{1})L_{2}K_{3}
\nonumber \\
K_{2}L_{2} &=&qL_{2}K_{2}-qQ_{1}^{-1}L_{1}^{\ast
}K_{1}+q^{-1}Q_{1}L_{1}K_{3})  \nonumber \\
K_{3}K_{3}^{\ast } &=&q^{-2}Q_{1}^{3}K_{3}^{\ast }K_{3}+L_{1}^{\ast }L_{1}
\nonumber \\
K_{3}L_{1} &=&q^{-1}Q_{1}^{2}L_{1}K_{3}  \nonumber \\
K_{3}L_{1}^{\ast } &=&qL_{1}^{\ast }K_{3}  \nonumber \\
K_{3}L_{2} &=&Q_{1}L_{2}K_{3}  \nonumber \\
L_{1}L_{1}^{\ast } &=&q^{2}Q_{1}^{-2}L_{1}^{\ast }L_{1}  \nonumber \\
L_{1}L_{2} &=&qQ_{1}^{-1}L_{2}L_{1}.  \nonumber
\end{eqnarray}

\noindent The Hopf algebra structure is given by the group Hopf algebra

\begin{equation}
\Delta (\mathbf{t})=\mathbf{t}\otimes \mathbf{t,\quad }\epsilon (\mathbf{t})=%
\mathbf{1,\quad }S(\mathbf{t)=t}^{-1}
\end{equation}

\noindent where the inverse matrix is given by

\begin{equation}
\mathbf{t}^{-1}{\small =\left(
\begin{array}{ccc}
L_{2}K_{1}^{\ast }-qQ_{1}^{-1}L_{1}^{\ast }K_{2}^{\ast } &
-Q_{1}^{-2}L_{2}K_{3}^{\ast }+q^{-1}Q_{1}^{-1}L_{1}K_{2}^{\ast } &
qQ_{1}^{-2}L_{1}^{\ast }K_{3}^{\ast }-q^{-1}L_{1}K_{1}^{\ast } \\
-Q_{1}^{2}L_{2}K_{3}+qQ_{1}L_{1}^{\ast }K_{2} &
L_{2}K_{1}-q^{-1}Q_{1}L_{1}K_{2} & q^{-1}Q_{1}^{2}L_{1}K_{3}-qL_{1}^{\ast
}K_{1} \\
q^{2}Q_{1}^{-1}K_{3}K_{2}^{\ast }-qK_{2}K_{1}^{\ast } & q^{-2}Q_{1}K_{3}^{%
\ast }K_{2}-q^{-1}K_{2}^{\ast }K_{1} & K_{1}^{\ast
}K_{1}-q^{-2}Q_{1}^{2}K_{3}^{\ast }K_{3}
\end{array}
\right) \delta ^{-1}.}
\end{equation}

\noindent The element $\delta $ which is defined to be

\begin{eqnarray}
\delta \equiv L_{2}K_{1}^{\ast }K_{1}-q^{-2}Q_{1}^{2}L_{2}K_{3}^{\ast
}K_{3}+L_{1}K_{3}K_{2}^{\ast }+L_{1}^{\ast }K_{3}^{\ast
}K_{2}-qQ_{1}^{-1}L_{1}^{\ast }K_{2}^{\ast
}K_{1}-q^{-1}Q_{1}L_{1}K_{2}K_{1}^{\ast }
\end{eqnarray}

\noindent has grouplike Hopf algebra structure

\begin{equation}
\Delta (\delta )=\delta \otimes \delta ,\quad \epsilon (\delta )=1,\quad
S(\delta )=\delta ^{-1}
\end{equation}

\noindent and satisfies

\begin{equation}
K_{1}\delta =\delta K_{1},\ K_{2}\delta =q^{-1}Q_{1}^{2}\delta K_{2},\
K_{3}\delta =q^{-2}Q_{1}^{4}\delta K_{3},\ L_{1}\delta =qQ_{1}^{-2}\delta
L_{1},\ L_{2}\delta =\delta L_{2}
\end{equation}

\noindent and their *-conjugates with $\delta ^{\ast }=\delta$ .

The braided Hopf algebra structure of the generalized oscillator (\ref
{eq:osc1}) implied by the quantum group covariance is given by the coproducts

\begin{equation}
\Delta (q^{N}) =q^{N}\otimes 1+1\otimes q^{N},\ \Delta (a)=a\otimes
1+1\otimes a, \ \Delta (a^{\ast }) =a^{\ast }\otimes 1+1\otimes a^{\ast },
\end{equation}

\noindent the counits

\begin{equation}
\epsilon(q^{N})= \epsilon(a)= \epsilon(a^{\ast })=0,
\end{equation}
\noindent

\noindent the antipodes

\begin{equation}
S(q^{N}) =-q^{N}, \ S(a) =-a, \ S(a^{\ast }) =-a^{\ast }
\end{equation}

\noindent and the braidings implied by (\ref{eq:rmprime2})

\begin{eqnarray}
\psi (q^{N}\otimes q^{N}) &=&q^{2}Q_{1}^{-1}q^{N}\otimes q^{N},  \nonumber
\label{eq:braid} \\
\psi (q^{N}\otimes a) &=&qQ_{1}^{-1}a\otimes q^{N},  \nonumber \\
\psi (a^{\ast }\otimes q^{N}) &=&qQ_{1}^{-1}q^{N}\otimes a^{\ast },
\nonumber \\
\psi (q^{N}\otimes a^{\ast }) &=&qa^{\ast }\otimes
q^{N}+Q_{1}^{-1}(q^{2}-Q_{1})q^{N}\otimes a^{\ast },  \nonumber \\
\psi (a\otimes q^{N}) &=&qq^{N}\otimes a+Q_{1}^{-1}(q^{2}-Q_{1})a\otimes
q^{N},  \nonumber \\
\psi (a\otimes a) &=&q^{2}Q_{1}^{-1}a\otimes a, \\
\psi (a^{\ast }\otimes a^{\ast }) &=&q^{2}Q_{1}^{-1}a^{\ast }\otimes a^{\ast
},  \nonumber \\
\psi (a\otimes a^{\ast }) &=&Q_{1}^{-1}(q^{2}-Q_{1})a\otimes a^{\ast
}+Q_{1}a^{\ast }\otimes a+q^{N}\otimes q^{N},  \nonumber \\
\psi (a^{\ast }\otimes a) &=&-q^{2}Q_{1}^{-2}q^{N}\otimes
q^{N}+q^{2}Q_{1}^{-2}a\otimes a^{\ast }.  \nonumber
\end{eqnarray}

\noindent In contrast to the three parameter deformed case where there is a
unique solution for the braidings, the two parameter deformed case $%
Q_{1}=q^{2}$ has three more solutions apart from the solution obtained by
substituting $Q_{1}=q^{2}$ into (\ref{eq:braid}). These solutions are

sol1:
\begin{eqnarray}
C_{1} &=&1,\ C_{2}=q^{2}\ ,C_{3}=0,\ C_{4}=0,\ C_{5}=q,\ C_{6}=0,\ C_{7}=0,
\nonumber  \label{eq:sol1} \\
C_{8} &=&q^{-2},\ C_{9}=0,\ C_{10}=q^{-1},\ C_{11}=0,\ C_{12}=0,\ C_{13}=0,\
C_{14}=-1
\end{eqnarray}

\bigskip

sol2:
\begin{eqnarray}
C_{1} &=&1,\ C_{2}=q^{2}\ ,C_{3}=0\ ,C_{4}=2\ ,C_{5}=q,\ C_{6}=0,\ C_{7}=0,
\nonumber \\
C_{8} &=&q^{-2},\ C_{9}=0,\ C_{10}=q^{-1},\ C_{11}=0,\ C_{12}=0,\ C_{13}=0,\
C_{14}=1
\end{eqnarray}

sol3:
\begin{eqnarray}
C_{1} &=&1,\ C_{2}=q^{2},\ C_{3}=0,\ C_{4}=0,\ C_{5}=q,\ C_{6}=0,\ C_{7}=0,
\nonumber  \label{eq:sol3} \\
C_{8} &=&q^{-2},\ C_{9}=-2q^{-2},\ C_{10}=q^{-1},\ C_{11}=0,\ C_{12}=0,\
C_{13}=0,\ C_{14}=1.
\end{eqnarray}

The $Q_{1}=q^{2}$ case is special not only because there are three more
solutions for the braidings, but also the $\mathbf{R}$-matrix is triangular (%
\textbf{\ }$\mathbf{R}_{12}^{-1}=\mathbf{R}_{21})$ and $S^{2}=id$ is
satisfied for the quantum group. We also note that in the general braided
Hopf algebra solutions given in Ref.4, only the solutions we give in this
section are related with the quantum group covariance.

We should also note that when $L_{1}=L_{1}^{\ast }=K_{2}=K_{2}^{\ast }=0$
and $L_{2}=1$ the transformation matrix is an element of $SU(1,1)\ $in the \
$q=Q_{1}=1$ limit.  Hence the group we define can be interpreted as
deformations of $SU(1,1).$

\section{ Subgroups and representations}

In the general form of the transformation of the generalized oscillator, the
invariance quantum group is a nine-parameter quantum group with three
deformation parameters. This quantum group has seven and five parameter
subgroups which we are going to discuss.

\textbf{A: } The seven parameter subgroup can be obtained by setting $%
L_{1}=L_{1}^{\ast }=0$ in (\ref{eq:relations1}). Then the consistency of the
relations requires $Q_{1}=q^{2}$ , i.e., for the oscillator

\begin{eqnarray}
aa^{\ast }-q^{2}a^{\ast }a &=&q^{2N}  \nonumber \\
aq^{N} &=&qq^{N}a \\
q^{N}a^{\ast } &=&qa^{\ast }q^{N}  \nonumber
\end{eqnarray}

\noindent the transformation

\begin{equation}
(
\begin{array}{ccc}
a & a^{\ast } & q^{N}
\end{array}
)^{\prime }=(
\begin{array}{ccc}
a & a^{\ast } & q^{N}
\end{array}
)\left(
\begin{array}{ccc}
K_{1} & K_{3}^{\ast } & 0 \\
K_{3} & K_{1}^{\ast } & 0 \\
K_{2} & K_{2}^{\ast } & L_{2}
\end{array}
\right)
\end{equation}

\noindent leaves the algebra covariant where the entries of the quantum
matrix satisfy

\begin{eqnarray}
K_{1}K_{1}^{\ast } &=&K_{1}^{\ast }K_{1},  \nonumber  \label{eq:sub1} \\
K_{1}K_{2} &=&q^{-1}K_{2}K_{1},  \nonumber \\
K_{1}K_{2}^{\ast } &=&qK_{2}^{\ast }K_{1},  \nonumber \\
K_{1}K_{3} &=&q^{-2}K_{3}K_{1},  \nonumber \\
K_{1}K_{3}^{\ast } &=&q^{2}K_{3}^{\ast }K_{1},  \nonumber \\
K_{1}L_{2} &=&L_{2}K_{1}, \\
K_{2}K_{2}^{\ast } &=&q^{2}K_{2}^{\ast }K_{2}+q^{2}K_{3}^{\ast
}K_{3}-K_{1}^{\ast }K_{1}+L_{2}^{2},  \nonumber \\
K_{2}K_{3} &=&q^{-1}K_{3}K_{2},  \nonumber \\
K_{2}K_{3}^{\ast } &=&q^{3}K_{3}^{\ast }K_{2},  \nonumber \\
K_{2}L_{2} &=&qL_{2}K_{2},  \nonumber \\
K_{3}K_{3}^{\ast } &=&q^{4}K_{3}^{\ast }K_{3},  \nonumber \\
K_{3}L_{2} &=&q^{2}L_{2}K_{3}.  \nonumber
\end{eqnarray}

\noindent The Hopf algebra structure is given by the group Hopf algebra,
i.e.,
\begin{equation}
\Delta (\mathbf{{t})={t}\otimes {t},\quad \epsilon ({t})=1,\quad S({t})={t}%
^{-1}}
\end{equation}

\noindent and the matrix inverse is

\begin{equation}
\mathbf{t}^{-1}=\left(
\begin{array}{ccc}
L_{2}K_{1}^{\ast } & -q^{-4}L_{2}K_{3}^{\ast } & 0 \\
-q^{4}L_{2}K_{3} & L_{2}K_{1} & 0 \\
K_{3}K_{2}^{\ast }-qK_{2}K_{1}^{\ast } & K_{3}^{\ast
}K_{2}-q^{-1}K_{2}^{\ast }K_{1} & K_{1}^{\ast }K_{1}-q^{2}K_{3}^{\ast }K_{3}
\end{array}
\right) \delta ^{-1}
\end{equation}
\noindent where the element

\begin{equation}
\delta \equiv L_{2}(K_{1}^{\ast }K_{1}-q^{2}K_{3}^{\ast }K_{3})
\end{equation}

\noindent is grouplike
\begin{equation}
\Delta (\delta )=\delta \otimes \delta ,\ \epsilon (\delta )=1,\ S(\delta
)=\delta ^{-1}
\end{equation}

\noindent and satisfies

\begin{equation}
K_{1}\delta =\delta K_{1},\ K_{2}\delta =q^{3}\delta K_{2},\ K_{3}\delta
=q^{6}\delta K_{3},\ L_{2}\delta =\delta L_{2},\ \delta ^{\ast }=\delta.
\end{equation}

To construct the representation, we take the generators of this algebra as
operators acting on some space. We first find the simultaneously
diagonalizible operators: the operators $L_{2}$, $K_{1}$ and $K_{1}^{\ast }$
commute among themselves and taking into account that $L_{2}^{\ast }=L_{2}$
we can take these operators as diagonal operators. We take the eigenvalue of
the Hermitian operator as
\begin{equation}
L_{2}\mid n\rangle =Aq^{n}\mid n\rangle
\end{equation}

\noindent where $A$ is a real constant. The relations of the algebra suggest
that

\begin{equation}
K_{2}\mid n\rangle \sim \mid n-1\rangle ,\quad K_{3}\mid n\rangle \sim \mid
n-2\rangle
\end{equation}

\noindent and hence we take the actions of the generators as
\begin{eqnarray}
K_{1}\mid n\rangle =k_{1,n}\mid n\rangle ,\quad K_{2}\mid n\rangle
&=&k_{2,n}\mid n-1\rangle \quad K_{3}\mid n\rangle =k_{3,n}\mid n-2\rangle \\
K_{1}^{\ast }\mid n\rangle =k_{1,n}^{\ast }\mid n\rangle ,\quad K_{2}^{\ast
}\mid n\rangle &=&k_{2,n+1}^{\ast }\mid n+1\rangle ,\quad K_{3}^{\ast }\mid
n\rangle =k_{3,n+2}^{\ast }\mid n+2\rangle .  \nonumber
\end{eqnarray}

\noindent Substituting these into (\ref{eq:sub1}) we obtain
\begin{equation}
k_{1,n}=Bq^{n},\quad k_{2,n}=Cq^{n},\quad k_{3,n}=Dq^{n}
\end{equation}

\noindent where B, C and D are complex constants satisfying
\begin{equation}
\left| B\right| ^{2}=A^{2}+q^{2}\left| D\right| ^{2}.
\end{equation}

\noindent We note that the representation is infinite dimensional.

\textbf{B:} The five parameter subgroup can be obtained by setting $%
L_{1}=L_{1}^{\ast }=K_{3}=K_{3}^{\ast }=0$ in the transformation (\ref
{eq:tfm1}), i.e., for the algebra

\begin{eqnarray}
aa^{\ast }-Q_{1}a^{\ast }a &=&q^{2N}  \nonumber \\
aq^{N} &=&qq^{N}a \\
q^{N}a^{\ast } &=&qa^{\ast }q^{N}  \nonumber
\end{eqnarray}

\noindent the transformation

\begin{equation}
(
\begin{array}{ccc}
a & a^{\ast } & q^{N}
\end{array}
)^{\prime }=(
\begin{array}{ccc}
a & a^{\ast } & q^{N}
\end{array}
)\left(
\begin{array}{ccc}
K_{1} & 0 & 0 \\
0 & K_{1}^{\ast } & 0 \\
K_{2} & K_{2}^{\ast } & L_{2}
\end{array}
\right)
\end{equation}

\noindent leaves the algebra covariant where the entries of the quantum
matrix satisfy

\begin{eqnarray}
K_{1}K_{1}^{\ast } &=&K_{1}^{\ast }K_{1},  \nonumber  \label{eq:sub2} \\
K_{1}K_{2} &=&qQ_{1}^{-1}K_{2}K_{1},  \nonumber \\
K_{1}K_{2}^{\ast } &=&q^{-1}Q_{1}K_{2}^{\ast }K_{1},  \nonumber \\
K_{1}L_{2} &=&L_{2}K_{1} \\
K_{2}K_{2}^{\ast } &=&Q_{1}K_{2}^{\ast }K_{2}-K_{1}^{\ast }K_{1}+L_{2}^{2}
\nonumber \\
K_{2}L_{2} &=&qL_{2}K_{2}.  \nonumber
\end{eqnarray}

\noindent The Hopf algebra structure is given by the group Hopf algebra
\begin{equation}
\Delta (\mathbf{{t})={t}\otimes {t},\quad \epsilon ({t})=1,\quad S({t})={t}%
^{-1}}.
\end{equation}
\noindent The inverse matrix is
\begin{equation}
\mathbf{{t}^{-1}}=\left(
\begin{array}{ccc}
L_{2}K_{1}^{\ast } & 0 & 0 \\
0 & L_{2}K_{1} & 0 \\
-qK_{2}K_{1}^{\ast } & -q^{-1}K_{2}^{\ast }K_{1} & K_{1}^{\ast }K_{1}
\end{array}
\right) \delta ^{-1}
\end{equation}

\noindent where

\begin{equation}
\delta \equiv L_{2}K_{1}^{\ast }K_{1}
\end{equation}

\noindent is grouplike
\begin{equation}
\Delta (\delta )=\delta \otimes \delta ,\quad \epsilon (\delta )=1,\quad
S(\delta )=\delta ^{-1}
\end{equation}

\noindent and satisfies
\begin{equation}
K_{1}\delta =\delta K_{1},\quad K_{2}\delta =q^{-1}Q_{1}^{2}\delta
K_{2},\quad L_{2}\delta =\delta L_{2}.
\end{equation}

Similar to the construction of the representation of the seven parameter
subgroup, the elements $L_{2}$ and $K_{1}$ can be taken as diagonal
operators and $K_{2}$ and $K_{2}^{\ast }$ can be taken as lowering and
raising operators respectively. We take the eigenvalue of the Hermitian
operator as
\begin{equation}
L_{2}\mid n\rangle =Aq^{n}\mid n\rangle
\end{equation}

\noindent and for the other operators we take

\begin{eqnarray}
K_{1}\mid n\rangle =k_{1,n}\mid n\rangle ,\ K_{2}\mid n\rangle =k_{2,n}\mid
n-1\rangle , \\
K_{1}^{\ast }\mid n\rangle =k_{1,n}^{\ast }\mid n\rangle ,\ K_{2}^{\ast
}\mid n\rangle =k_{2,n+1}^{\ast }\mid n+1\rangle .  \nonumber
\end{eqnarray}

\noindent Substituting these into (\ref{eq:sub2}) we obtain
\begin{equation}
k_{1,n}=B\QDOVERD( ) {Q_{1}}{q}^{n},\quad \left| k_{2,n}\right| ^{2}=A^{2}%
\dfrac{Q_{1}^{n}-q^{2n}}{Q_{1}-q^{2}}-\left| B\right| ^{2}\dfrac{%
Q_{1}^{n}-\QDOVERD( ) {Q_{1}}{q}^{2n}}{Q_{1}-\QDOVERD( ) {Q_{1}}{q}^{2}}
\end{equation}

\noindent where $A$ is real and $B$ is complex. The quadratic Casimir of the
algebra which is found to be
\begin{equation}
C=K_{1}^{\ast }K_{1}+(q^{-2}-1)K_{2}^{\ast }K_{2}+q^{-2}L_{2}^{2}
\end{equation}
\noindent has the eigenvalue
\begin{equation}
C\mid n\rangle =(A^{2}+q^{-2}\left| B\right| ^{2})\mid n\rangle .
\end{equation}

In the algebra (\ref{eq:sub2}) when we identify

\begin{equation}
L_{2}\equiv q^{H},\ \ K_{1}=K_{1}^{\ast }\equiv q^{-H},\ \ K_{2}\equiv
(q-q^{-1})^{1/2}X_{-},\ \ K_{2}^{\ast }\equiv (q-q^{-1})^{1/2}X_{+},\ \
Q_{1}=1
\end{equation}

\noindent the algebra turns out to be

\begin{equation}
q^{H}X_{\pm }=q^{\pm 1}X_{\pm }q^{H}, \quad X_{+}X_{-}-X_{-}X_{+}= \dfrac{%
q^{2H}-q^{-2H}}{q-q^{-1}}
\end{equation}

\noindent which generates $U_{q}(su(2))$ with the *-structure

\begin{equation}
(q^{H})^{\ast }=q^{H},\quad (X_{\pm })^{\ast }=X_{\mp }.
\end{equation}

\noindent The transformation

\begin{equation}
(
\begin{array}{ccc}
a & a^{\ast } & q^{N}
\end{array}
)^{\prime }=(
\begin{array}{ccc}
a & a^{\ast } & q^{N}
\end{array}
)\left(
\begin{array}{ccc}
q^{-H} & 0 & 0 \\
0 & q^{-H} & 0 \\
X_{-} & X_{+} & q^{H}
\end{array}
\right)
\end{equation}

\noindent leaves the commutation relations between $a,a^{\ast }$ and $q^{N}$
invariant. The matrix multiplication gives the comultiplication of the
generators of the algebra

\begin{equation}
\Delta (q^{H})=q^{H}\otimes q^{H},\quad \Delta (q^{-H})=q^{-H}\otimes
q^{-H},\quad \Delta (X_{\pm })=X_{\pm }\otimes q^{-H}+q^{H}\otimes X_{\pm }
\end{equation}

\noindent and the counit map

\begin{equation}
\epsilon (q^{\pm H})=1,\quad \epsilon (X_{\pm })=0
\end{equation}

\noindent is a map to the identity of the quantum group. The antipodes

\begin{equation}
S(q^{\pm H})=q^{\mp H},\quad S(X_{\pm })=-q^{\mp 1}X_{\pm }
\end{equation}

\noindent give the inverse of the transformation matrix.

\section{CONCLUSION}

At the level of a single oscillator, the deformation of the oscillator
results in a noncommutativity in the algebra whose coaction leaves the
oscillator algebra covariant. At the level of two (or more) oscillators,
this induces a noncommutativity (called outer noncommutativity) between
independent copies. This noncommutativity is described by the braiding
relations. The discussion of the $n$-fold braided tensor product for $q$%
-Heisenberg algebra is done in Ref.5. Hence we give not only a
generalization of the oscillator algebra but also the interaction pattern of
these oscillators among themselves via the quantum group covariance. This
may contribute to understand the possible connections between quantum groups
and nonextensive statistical mechanics$^{6}$. In this work we use real
deformation parameters, however, the fractional supersymmetric structures
require at least one of the deformation parameters to be a root of unity$^{7}
$ as the generalization of the (-1) factor in the fermionic case. One more
thing which deserves a separate study is the decoupling of the oscillators
or the unbraiding transformations$^{8}$.

The braided covector algebras covariant under quantum groups are also
covariant under the braided groups obtained from quantum groups by a
transmutation process$^{3}$. The main ingredient of this construction is the
$\mathbf{R}$-matrix. Hence the $\mathbf{R}$-matrix we found defines a new
braided group which we do not consider here.

\pagebreak

\textbf{References}

\begin{description}
\item  {$^{1}$}R. Jagannathan, R. Sridhar, R. Vasudevan, S. Chaturvedi, M.
Krishnakumari, P. Shanta, and V. Srinivasan, J. Phys. A: Math. Gen. \textbf{%
25}, 6429 (1992); J. Bertrand, M. Iracastaud, J. Phys. A: Math. Gen. \textbf{%
30}, 2021 (1997); M. Arik , AS Arikan, J. Math. Phys. \textbf{42}, 2388
(2001).

\item  {$^{2}$}S. Majid, Algebras and Hopf algebras in braided categories,
Lecture Notes in Pure and Applied Mathematics, vol 158 (New York: Marcel
Dekker), pp. 55-105 (1994).

\item  {$^{3}$}S.Majid, \textit{Foundations of Quantum Group Theory}
(Cambridge: Cambridge University Press), (1995) S. Majid, Beyond
supersymmetry and quantum symmetry (an introduction to braided groups and
braided matrices) \textit{Quantum Groups, Integrable Statistical Models and
Knot Theory} edited by M. L. Ge and H. J. Vega (Singapore: World
Scientific), pp.231-282 (1993).

\item  {$^{4}$}A. Yildiz, J. Math. Phys. \textbf{43}, 1668 (2002).

\item  {$^{5}$}W. K. Baskerville and S. Majid, J. Math. Phys. \textbf{34},
3588 (1993).

\item  {$^{6}$}M.R. Ubriaco, Phys. Lett. A \textbf{283}, 157 (2001); S. Abe,
\textit{ibid} \textbf{244}, 229 (1998);M. Arik, J. Kornfilt, A. Yildiz,%
\textit{ibid} \textbf{235}, 318 (1997); C. Tsallis, \textit{ibid} \textbf{195%
}, 329 (1994).

\item  {$^{7}$}R.S. Dunne, A.J. Macfarlane, J.A. deAzcarraga, J.C.P. Bueno,
Int. J. Mod. Phys. A, \textbf{12}, 3275 (1997); R.S. Dunne, J. Math. Phys.
\textbf{40}, 1180 (1999); H. Ahmedov and O. F. Dayi, J. Phys. A \textbf{32},
6247 (1999); H. Ahmedov, A. Yildiz and Y. Ucan, J. Phys. A \textbf{34}, 6413
(2001).

\item  {$^{8}$}G. Fiore, H. Steinacker and J. Wess, Mod. Phys. Lett. A
\textbf{16}, 261 (2001).
\end{description}

\end{document}